\documentclass[a4paper,11pt]{article}
\pdfoutput=1 

\usepackage{jheppub} 

\usepackage[T1]{fontenc} 
\usepackage{subcaption}
\usepackage{amsmath}
\usepackage{amssymb}
\usepackage{optidef}
\usepackage{doi}

\title{\boldmath A Note on Ribbon-based Biharmonic Surface Patches}

\author[a]{Márton Vaitkus}

\affiliation[a]{Budapest University of Technology and Economics,\\2 Magyar Tud\'{o}sok K\"{o}r\'{u}tja, Budapest, Hungary}

\emailAdd{vaitkus@iit.bme.hu}

\toccontinuoustrue
\makeatletter
\gdef\@fpheader{}
\makeatother

\abstract{In this short note we describe a simple adaptation of biharmonic surfaces to interpolate boundary cross-derivatives given in ribbon form, and compare with the recently proposed Generalized B-spline patches.}

\begin{document} 
\maketitle
\flushbottom

\section{Introduction}

Creating multi-sided general-topology surface patches that connect to their neighbourhood with geometric continuity \citep{Kiciak:2016} is a widely studied problem in computer-aided geometric design (CAGD). A large variety of methods for creating multi-sided surfaces have been proposed, see e.g. \citep{Peters:2019} and \citep{Vaitkus:2021} for overviews. While some progress has been made regarding mutlti-sided implicit surfaces \citep{Sipos:2020}, in practice parametric representations are generally preferred. Recently developed methods \citep{Varady:2020,Vaitkus:2021,Martin:2022} even allow for surfaces to be defined over domains that are bounded by arbitrary curves in the plane and are possibly multiply connected. Such surfaces are generally defined using boundary \emph{ribbons}, i.e. bi-parametric surfaces defining positional and cross-derivative information to be interpolated along the boundaries. For transfinite interpolation methods \citep{Sabin:1996,Dyken:2009,Varady:2011,Salvi:2014} ribbons can be defined as general (possibly procedurally generated) bi-parametric surfaces. Other methods are based on control points \citep{Varady:1991,Varady:2016,Varady:2020,Hettinga:2018,Hettinga:2020,Vaitkus:2021,Martin:2022}, in which case ribbons must be piecewise-rational surfaces in Bézier or B-spline form. As an example of the state-of-the-art, Generalized B-Spline (GBS) patches \citep{Vaitkus:2021} employ local re-parameterizations to map a general curved domain onto the domain of the B-spline ribbons, and pull back blending functions naturally associated with the ribbons onto the common domain, and then -- after applying suitable corrections to enforce interpolation of boundary derivatives and partition of unity -- each point of the GBS patch can be evaluated as a convex combination of the ribbon control points.

As an alternative to such direct methods, surfaces could also be defined  as solutions of functional optimization problems, or partial differential equations that express some measure of fairness or quality for the shape (such as integrated curvature or curvature variation). Solving nonlinear optimization problems or PDEs is generally only possible via computationally expensive iterative methods \citep{Moreton:1992,Schneider:2000,Pan:2015,Soliman:2021}. A more efficient alternative is to use quadratic energies, or linear PDEs which lead to simple linear systems of equations to solve \citep{Botsch:2004,Jacobson:2010,Andrews:2011}. Generally, such linear variational surfaces \citep{Botsch:2007} are defined in terms of normal derivatives along their boundary. In this short paper, we describe how to relate biharmonic surfaces as proposed by \citep{Jacobson:2010} to ribbon-based approaches and compare the resulting generalized Hermite interpolants to Generalized B-spline surfaces.

\section{Biharmonic interpolation and its mixed discretization}
In this chapter, we follow the derivations of \citep{Jacobson:2010}. Over a 2D domain $\Omega \subset \mathbb{R}^2$, a generalization of (cubic) Hermite interpolation called \emph{biharmonic} interpolation can be formulated as the constrained optimization of the \emph{Laplacian energy} (closely related -- but not always identical -- to the so-called thin-plate energy \citep{Stein:2018}):

\begin{mini}
    	{u}{\int_{\Omega} \left(\Delta u(x,y)\right)^2 dx dy}{}{}
    	\addConstraint{u(x,y) }{= u_0(x,y),\ (x,y) \in \partial \Omega}
            \addConstraint{\frac{\partial u}{\partial \mathbf{n}} }{= d_0(x,y)\ (x,y) \in \partial \Omega}
\end{mini}
where $\mathbf{n}(x,y)$ is the inwards normal at the boundary. We choose to discretize this problem using (mixed) finite elements, but note that alternatives such as boundary element \citep{Weber:2012} and Monte Carlo \citep{Sawhney:2020} methods could also be used to compute an approximate solution to the biharmonic problem.

Proceeding with a mixed formulation and introducing the auxiliary variable $v = \Delta u$ and Lagrange multiplier $\lambda$, we get the mixed Lagrangian:
\begin{align}
\min\limits_{u,v} \max\limits_{\lambda} \int_{\Omega} v^2 + \lambda\left(v - \Delta u\right) dx dy.
\end{align}
Applying Green's formula we derive the following saddle point problem:
\begin{align}
\min\limits_{u,v} \max\limits_{\lambda} \int_{\Omega} v^2 dx dy + \int_{\Omega} \lambda v\ dx dy + \int_{\Omega} \nabla \lambda \cdot \nabla u\ dx dy - \int_{\partial \Omega} \lambda \frac{\partial u}{\partial \mathbf{n}}\ dx dy.
\end{align}

Discretizing $\Omega$ with a triangulation, and approximating the unknowns as piecewise-linear functions 
\begin{align}
u(x) &= \sum_{i=1}^{n}u_i\varphi_i(x), \\
v(x) &= \sum_{i=1}^{n}v_i\varphi_i(x), \\
\lambda(x) &= \sum_{i=1}^{n}\lambda_i\varphi_i(x),
\end{align}
and differentiating with respect to the unknown coefficients we get a linear system 
\begin{align} 
\begin{bmatrix}
-\mathbf{M} & \mathbf{L} \\
\mathbf{L} & \mathbf{0}
\end{bmatrix} 
\begin{bmatrix}
\mathbf{v} \\
\mathbf{u}
\end{bmatrix} 
=
\begin{bmatrix}
\mathbf{b} \\
\mathbf{0}
\end{bmatrix} 
\end{align}
where 
\begin{align} 
\mathbf{L}_{ij} = \int_\Omega \nabla \varphi_i \cdot \nabla \varphi_j\ dx dy 
\end{align}
is the well-known cotangent Laplace matrix, 
\begin{align}
    \mathbf{M}_{ij} = \int_\Omega \varphi_i \varphi_j\ dx dy
\end{align}  
is the mass matrix, which is often approximated by a diagonal ("lumped mass") matrix containing the sum of each row (equal to the Voronoi area around each vertex). The right hand side of the equation is defined as
$\mathbf{b} = -\mathbf{L}_{\Omega,\partial \Omega} \mathbf{u}_0 - \mathbf{N}\mathbf{d}_0 $
where $\mathbf{L}_{\Omega,\partial \Omega}$ is the Laplacian matrix restricted to the boundary, $\mathbf{N}_{ij} = \int_{\partial \Omega} \varphi_i \varphi_j\ dx dy $ is the boundary "area" (length) associated to a given vertex, and $\mathbf{u}_0,\mathbf{d}_0$ are sampled values of the boundary conditions.

The reduced form 
\begin{align}   
 \left( \mathbf{L}\mathbf{M}^{-1}\mathbf{L} \right)\mathbf{u} = \mathbf{L}\mathbf{M}^{-1}\mathbf{b} 
\end{align}
is similar to an alternative discretization of the biharmonic equation  \citep{Botsch:2004}. The vertex values can be computed as a linear combination of the boundary values and derivatives:
\begin{align}
    \mathbf{u} = \left( \mathbf{L}\mathbf{M}^{-1}\mathbf{L} \right)^{-1} \mathbf{L}\mathbf{M}^{-1}\mathbf{b} = \mathbf{H}_0 \mathbf{u}_0 + \mathbf{H}_1 \mathbf{d_0},
\end{align}
where the columns of the matrices $\mathbf{H}_0,\mathbf{H}_1$ can be interpeted as generalized cubic Hermite blending functions.

\section{Ribbon-based biharmonic patches}

In general, we have boundary conditions prescribed in the form of procedural or tensor-product Bernstein-Bézier or B-spline ribbons associated with our curved domain $\Omega$:
\begin{align}
\mathbf{R}_i:\  & [0,1]^2 \rightarrow\mathbb{R}^3 \\
& (s,h) \mapsto (x,y,z)
\end{align}
In particular, we want to use $\mathbf{R}(s,0)$ and $\frac{\partial \mathbf{R}} {\partial h}$ as our boundary conditions for the biharmonic interpolation. We locally (re)parameterize each ribbon onto the curved domain $\Omega$, using e.g. harmonic interpolation \citep{Joshi:2007,Vaitkus:2021}
\begin{align}
r_i:\  &\Omega \rightarrow [0,1]^2 \\
&(u,v) \mapsto (s_i,h_i)
\end{align}
As a consequence, the normal derivatives must be computed by also taking into account the gradient of the local parameterizations:
\begin{align}
   \frac{\partial\mathbf{R}_i}{\partial\mathbf{n}} = \frac{\partial s_i}{\partial \mathbf{n}}\frac{\partial \mathbf{R}_i}{\partial s} + \frac{\partial h_i}{\partial \mathbf{n}}\frac{\partial \mathbf{R}_i}{\partial h}
\end{align}
with e.g. $ \frac{\partial s_i}{\partial \mathbf{n}} = \nabla s_i \cdot \mathbf{n}$, where $\nabla s_i = \begin{bmatrix} \frac{\partial s_i}{\partial u} & \frac{\partial s_i}{\partial v} \end{bmatrix}$ can be approximated over the domain mesh with any of the usual methods \citep{Mancinelli:2019}. Note that unlike other ribbon-based patches that directly use the parameterization in their definition \citep{Salvi:2014,Vaitkus:2021}, for biharmonic patches only the boundary derivatives of the reparameterization are utilized, similarly to the usual sufficient conditions for geometric continuity \citep{Kiciak:2016}.

\section{Discussion}
\begin{figure}[h]
	\centering
	\begin{subfigure}{0.45\textwidth}
		\includegraphics[width=\textwidth,keepaspectratio]{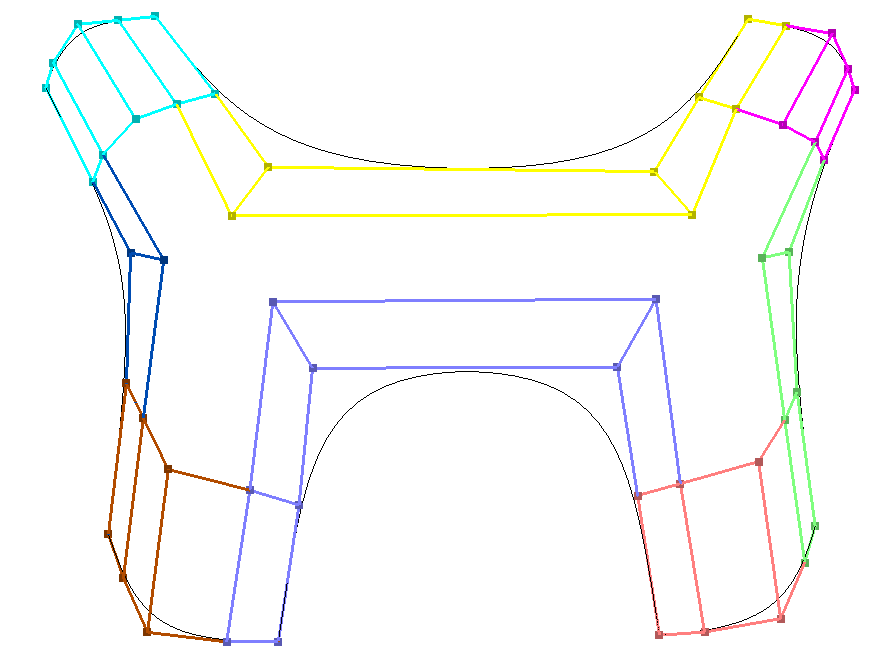}
		\subcaption{Input ribbons}
	\end{subfigure}

 	\begin{subfigure}{0.45\textwidth}
		\includegraphics[width=\textwidth,keepaspectratio]{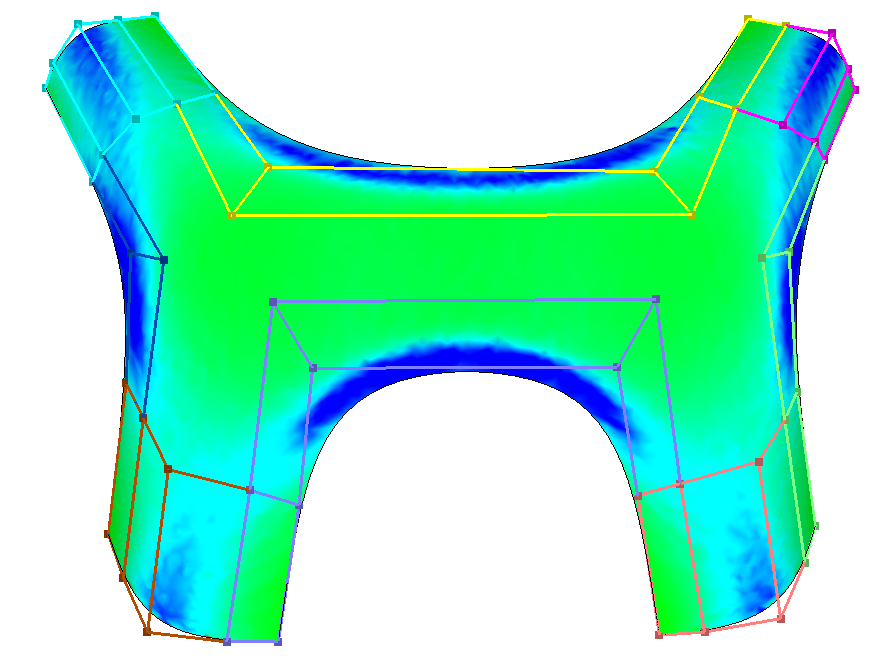}
		\subcaption{GBS patch}
	\end{subfigure} \
 	\begin{subfigure}{0.45\textwidth}
		\includegraphics[width=\textwidth,keepaspectratio]{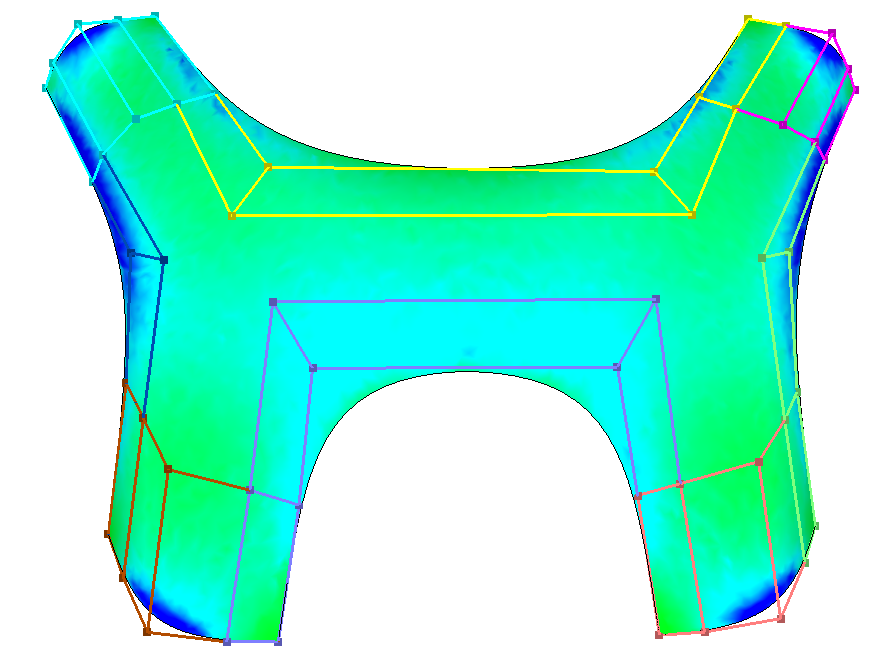}
		\subcaption{Biharmonic patch}
	\end{subfigure}
	\caption{Comparison of vertex blend surfaces. Colors indicate mean curvature.}
	\label{fig:vb}
\end{figure}

\begin{figure}[h]
	\centering
	\begin{subfigure}{0.45\textwidth}
		\includegraphics[width=\textwidth,keepaspectratio]{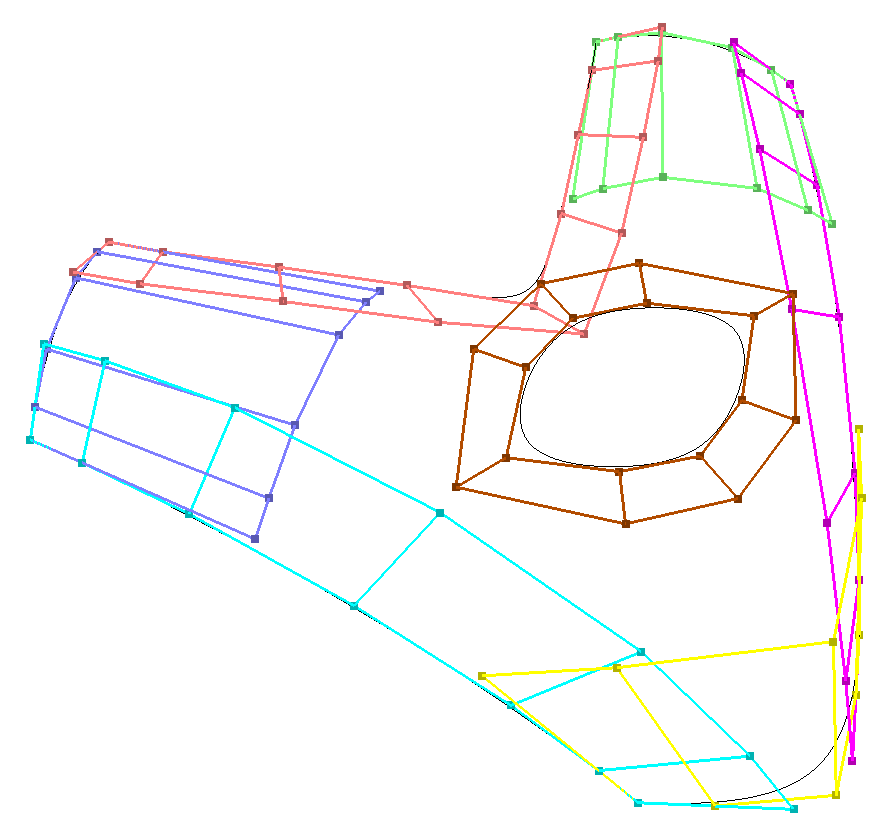}
		\subcaption{Input ribbons}
	\end{subfigure}

 	\begin{subfigure}{0.45\textwidth}
		\includegraphics[width=\textwidth,keepaspectratio]{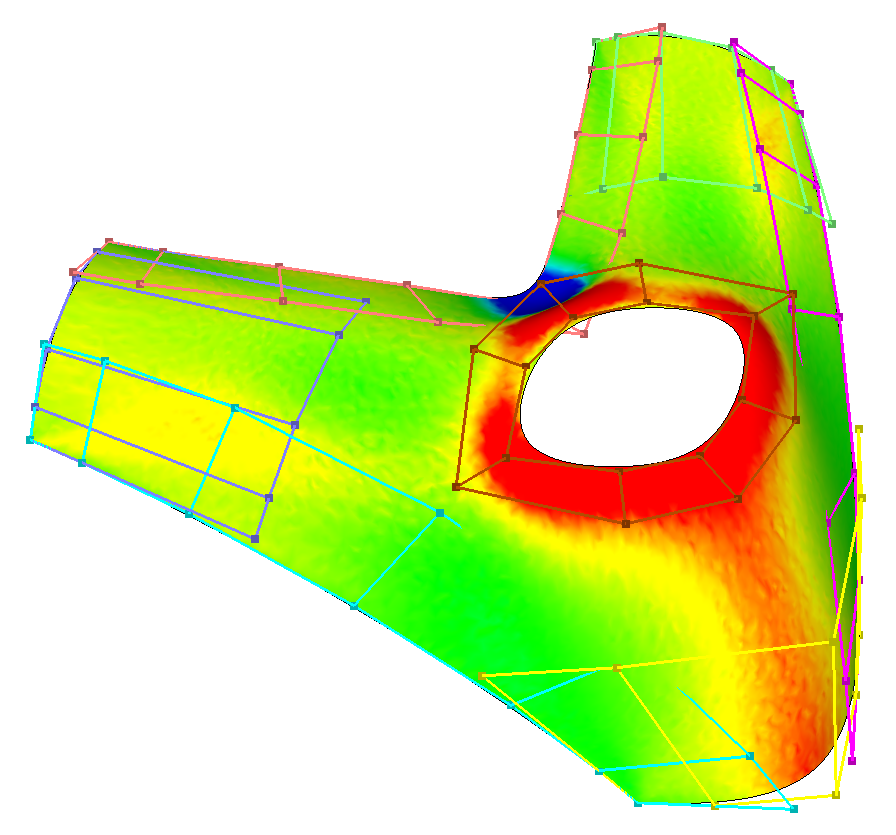}
		\subcaption{GBS patch}
	\end{subfigure} \
 	\begin{subfigure}{0.45\textwidth}
		\includegraphics[width=\textwidth,keepaspectratio]{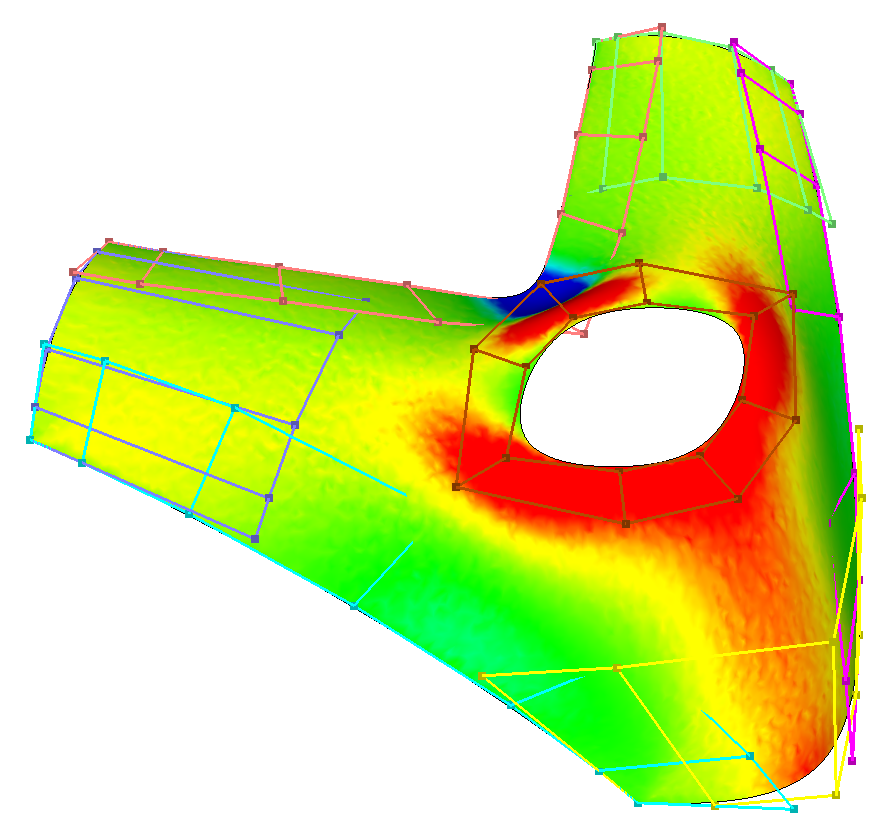}
		\subcaption{Biharmonic patch}
	\end{subfigure}
	\caption{Comparison of multiply connected surfaces. Colors indicate mean curvature.}
	\label{fig:ml}
\end{figure}

We compare B-spline ribbon-based biharmonic patches with Generalized B-spline surfaces of \citep{Vaitkus:2021} in \autoref{fig:vb} and \autoref{fig:ml}. As expected, biharmonic patches usually (although not always) have a more even curvature distribution than GBS patches, due to their energy-based formulation. However, our experience is that biharmonic patches are less responsive to control point movements than more explicit methods such as GBS patches (likely due to the built-in fairness), which might make make it difficult to do fine-grained manual design. This supports earlier claims by \citep{Peters:2003,Peters:2008}.

\begin{figure}[h]
	\centering
	\begin{subfigure}{0.45\textwidth}
		\includegraphics[width=\textwidth,keepaspectratio]{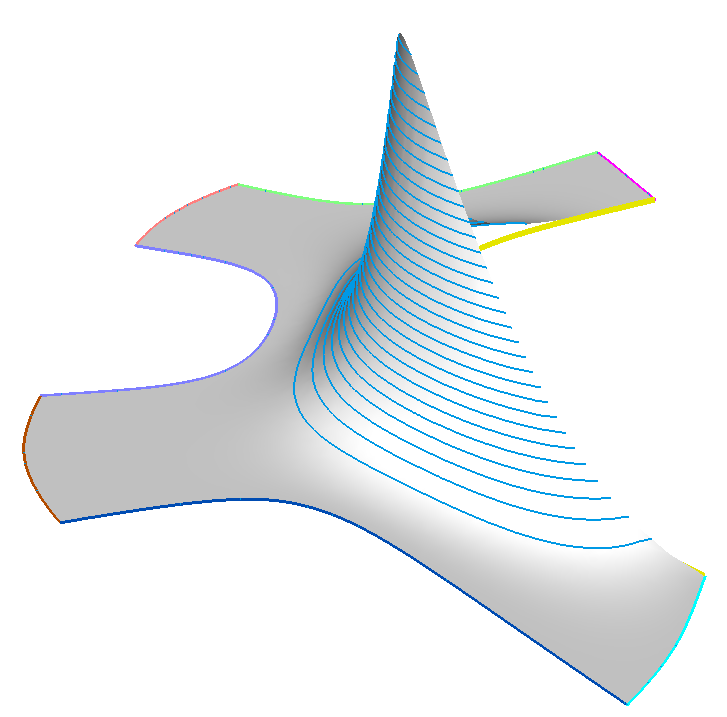}
		\subcaption{GBS (boundary CP)}
	\end{subfigure}\ 
 	\begin{subfigure}{0.45\textwidth}
		\includegraphics[width=\textwidth,keepaspectratio]{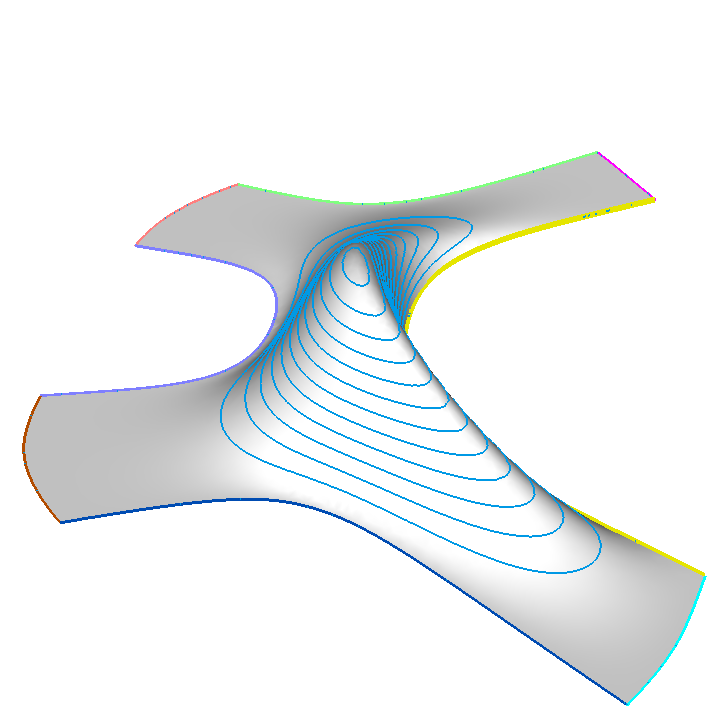}
		\subcaption{GBS (tangent CP)}
	\end{subfigure} 
 
	\begin{subfigure}{0.45\textwidth}
		\includegraphics[width=\textwidth,keepaspectratio]{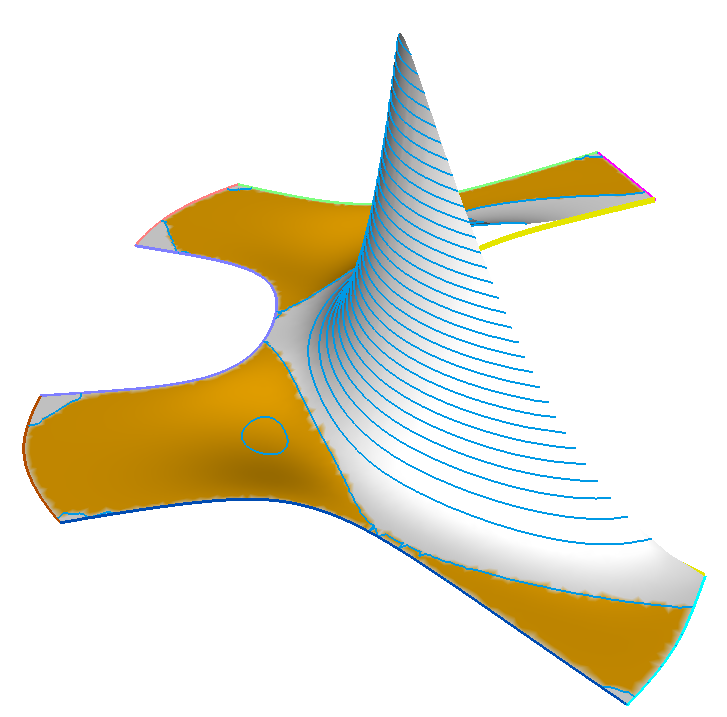}
		\subcaption{Biharmonic (boundary CP)}
	\end{subfigure}\ 
 	\begin{subfigure}{0.45\textwidth}
		\includegraphics[width=\textwidth,keepaspectratio]{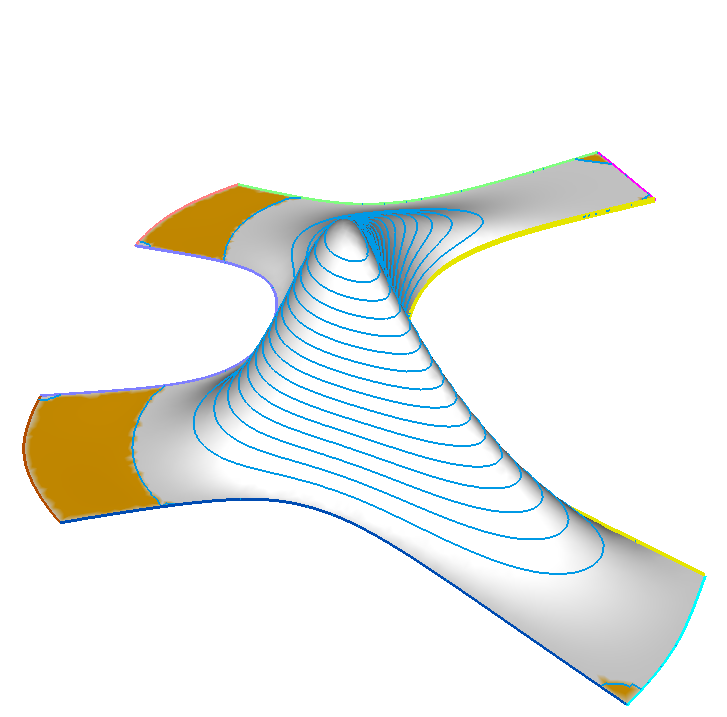}
		\subcaption{Biharmonic (tangent CP)}
	\end{subfigure} 
	\caption{Comparison of blend function level surfaces. Yellow color indicates negative values.}
	\label{fig:blend}
\end{figure}

For biharmonic patches, blending functions can be pre-evaluated by substituting Kronecker deltas for the positions and normal derivatives and converting from a Hermite to a Bernstein representation of cross-derivatives. The resulting blending functions are compared with those of GBS patches on \autoref{fig:blend}. As can be seen, the functions are often fairly similar, although biharmonic blending functions can be negative, while GBS blending functions are guaranteed to form a convex partition of unity.

Instead of the presented control point-based approach, one could also interpret ribbons of biharmonic surfaces in a Hermite-like manner as in e.g. \citep{Varady:2011,Salvi:2014}. We mention that ribbon-based triharmonic surfaces could also be derived following \citep{Jacobson:2010}, based on the second derivatives of the reparameterization -- although the accurate calculation of these over triangulated domains might be challenging. 

\bibliographystyle{JHEP}
\bibliography{biharmonic}
\end{document}